\documentclass[11pt]{amsart}
\usepackage{graphicx}
\usepackage{color}
\usepackage[english]{babel}
\usepackage[latin1]{inputenc}
\usepackage{amsmath}
\usepackage{amssymb}
\newtheorem{teo}{Theorem}[section]
\newtheorem{cor}[teo]{Corollary}
\newtheorem{lema}[teo]{Lemma}
\newtheorem{prop}[teo]{Proposition}
\theoremstyle{definition}
\newtheorem{defn}[teo]{Definition}
\theoremstyle{remark}
\newtheorem{rem}[teo]{Remark}
\newtheorem{ex}[teo]{Example}
\numberwithin{equation}{section}
\newcommand{\N}{\mathbb{N}}

\newcommand{\dem}{\noindent {\textsf{Proof.} }}

\author[Miralles]{Alejandro Miralles$^1$}\address{Alejandro Miralles.
Instituto Universitario de Matemáticas y Aplicaciones de Castellón (IMAC), Universitat Jaume I de Castelló (UJI), Castelló, Spain. \emph{e}.mail: mirallea@uji.es}

\author[Torres]{Damià Torres}
\address{Damià Torres. Universitat Politècnica de Catalunya - BarcelonaTech (UPC), Barcelona, Spain. 
\emph{e}.mail: damia.torres@estudiant.upc.edu}

\thanks{$^1$ Partially supported by MINECO project MTM2014-53241, Generalitat Valenciana project AICO/2016/030 and Universitat Jaume I project P1-1B2014-35 }
\begin{document}

\begin{abstract}
We introduce the sequence $(a_n) \subset (0,1]$ and prove that the asymptotic behaviour of $\sum_{k=1}^n a_k$ is the same than $\pi(n)$, the prime-counting function. We also obtain that $\pi(n) \sim n a_n$ and  we estimate $\frac{1}{a_n}-\frac{n}{\pi(n)}$ showing that $\lim_{n \rightarrow \infty} \frac{1}{a_n}-\frac{n}{\pi(n)}$ is convergent.
\end{abstract}

\title[Counting primes by sums of frequencies]{Counting primes by sums of frequencies}

\maketitle

\section{Introduction and background}
\subsection{Introduction}

In this work we introduce the numbers $a_n \in \left( 0,1\right]$ for any $n \in \mathbb{N}$ in two equivalent ways: first we use a sieve method which results in a subset $A_n \subset \{1,2,\cdots n\}$ and then we consider the frequency $a_n=\frac{|A_n|}{n}$. We prove that this is equivalent to consider the recurrence sequence $a_n$ given by $a_1=1$ and 
$$a_{n}=a_{n-1}\left(1-\frac{a_{n-1}}{n}\right)$$
\noindent as in Definition \ref{def pr}. We will prove that the asymptotic behaviour of the partial sums $c(n)=\sum_{k=1}^n a_n$ is the same than the prime counting function given by 
$$\pi(n)=\#\{p \in \mathbb{N} : p \mbox{ is prime and } p \leq n \}$$

\noindent for any $n \in \mathbb{N}$. We will also conclude other related results related to $\pi(n)$ as now that $\pi(n) \sim n a_n$ and  we will estimate $\frac{1}{a_n}-\frac{n}{\pi(n)}$ showing that $\lim_{n \rightarrow \infty} \frac{1}{a_n}-\frac{n}{\pi(n)}$ is convergent to a real number between $2$ and $3$.

\subsection{Classical results.} How prime numbers are distributed among positive integers plays an important role in number theory. Let $\pi(x)$ be the prime-counting function that gives the number of primes less than or equal to $x$ for any real number x.

The Prime Number Theorem gives the asymptotic behaviour of distribution of prime numbers. It was proved by Hadamard and de la Vallée-Poussin in 1896 using complex theory (see \cite{H} and \cite{V}). A more elementary proof was given by Erdös and Selberg in 1948 (see \cite{E} and \cite{S}). This classical result states that
$$\pi(x) \sim \frac{x}{\log x}$$

\noindent and also that
$$\pi(x) \sim Li(x),$$
\noindent where the asymptotic notation $f(x) \sim g(x)$ means 
$$\lim_{x \rightarrow \infty} \frac{f(x)}{g(x)}=1,$$

\noindent where $x$ could denote either a real or a positive integer and the Logarithmic Integral Function $Li(x)$ is given by
$$Li(x)= \int_{2}^x \frac{dt}{\log t} dt,$$
\noindent for any $x \geq 2$. We will consider the Harmonic numbers $H_n$ given by
$$H_n=1+\frac{1}{2}+\cdots+\frac{1}{n}$$

\noindent for any $n \in \mathbb{N}$. It is well-known that 
\begin{eqnarray} \label{desig2}
0< H_n -\log n - \gamma < \frac{1}{n},
\end{eqnarray}
\noindent where $\gamma \approx 0'577$ denotes the Euler-Mascheroni constant. \medskip

\subsection{Approaching the frequency $a_n$: sieving the positive integers}
Fix $n$ in $\N$, consider a \textit{good integer} $M$ (to be determined later) and let $A_0=\{1,2,\cdots,M\}$. For any $1 \leq k \leq n$, we will sieve some numbers up to $M$ and we will obtain a sequence of decreasing subsets $A_0 \supseteq A_1 \supseteq A_2 \supseteq \cdots \supseteq A_n$ and a decreasing sequence of real numbers $a_1 \geq a_2 \geq \cdots \geq a_n$ given by the frequency $a_k=\frac{|A_k|}{M}$ for any $k$. 

\begin{ex}
Fix $n=4$. We will study the frequencies $a_k$ for $k=1,2,3,4$. Consider $M=576$ and $A_0=\{1,2, \cdots, 576 \}$. \medskip

\textbf{Step 1.} For $k=1$ we start from the set $A_0$ and choose the subset $A_1$ of all the multiples of $1$, that is, the whole $A_0$. The frequency $a_1$ is given by $a_1=\frac{|A_1|}{|A_0|}=1$. \medskip

\textbf{Step 2.} For $k=2$, we start from the previous set $A_1$ of $|A_1|=576$ integers. We consider the subset $B_2$ by keeping every number in $A_1$ from the first one by counting up in increments of $2$ and crossing out the remaining ones, that is, $B_2=\{1,3,5, \cdots,575\}$, so $|B_2|=288$. Now we apply the previous step (step 1) to $B_2$, that is, we consider all the elements of $B_2$, which gives the set $C_2=\{1,3,5,\cdots,575\}$. The set $A_2$ is given by $A_1 \setminus C_2$ and, since $C_2 \subset A_1$, we have that $|A_2|=288$ and the frequency $a_2$ is given by $\frac{|A_2|}{M}=288/576=0.5$. \medskip

\textbf{Step 3.} For $k=3$, we start from the previous set $A_2$ of $|A_2|=288$ integers, that is, $A_2=\{2, 4, 6,  \cdots, 576\}$. We consider the subset $B_3$ by keeping every number in $A_2$ from the first one by counting up in increments of $3$ and crossing out the remaining ones, that is, $B_3=\{2,8,14,20,\cdots, 566, 572\}$, so $|B_3|=288/3=96$. We apply the previous step (step 2) to $B_3$, that is, we keep numbers in $B_3$ from the first one by counting up in increments of $2$ and crossing up the remaining ones, that is, $C_3=\{2,14,26, \cdots,  566\}$ satisfying $|C_3|=\frac{|A_2|}{M}|B_3|=48$ numbers. The set $A_3$ is given by $A_2 \setminus C_3$ and, since $C_3 \subset A_2$, we have that $|A_3|=288-48=240$ and the frequency $a_3$ is given by $\frac{|A_3|}{M}=240/512=5/12$. \medskip

\textbf{Step 4.} For $n = 4$, we start from the previous set $A_3$ of $|A_3|=240$ integers. We consider the subset $B_4$ by keeping every number in $A_3$ from the first one by counting up in increments of $4$ and crossing out the remaining ones, so $|B_4|=60$. We apply the previous step (step 3) to $B_4$ which gives the set $C_4$ satisfying $|C_4|=\frac{|A_3|}{M} |B_4|=5/12\cdot60=25$ numbers. The set $A_4$ is given by $A_3 \setminus C_4$ and since $C_4 \subset A_3$, we have that $|A_4|=240-25=215$ and the frequency $a_4$ is given by $\frac{|A_4|}{M}=215/576$. 
\end{ex}

To define the frequency $a_k$ for any $k \geq 4$, fix $n \geq k$ and consider a good integer $M$ for $n$ (to be determined later). Suppose $A_{k-1}$ and $a_{k-1}=\frac{|A_{k-1}|}{M}$ are defined. \medskip

\textbf{Step $k$.} We start from the previous set $A_{k-1}$ of $|A_{k-1}|=a_{k-1} M$ integers. We consider the subset $B_k$ by keeping every number in $A_{k-1}$ from the first one by counting up in increments of $k$ and crossing out the remaining ones. We apply the previous step (step k-1) to $B_k$ which gives the set $C_k$ satisfying $|C_k|=\frac{|A_{k-1}|}{M}|B_k|=\frac{|A_{k-1}|^2}{M k}=a_{k-1} |B_k|$. The set $A_k$ is given by $A_{k-1} \setminus C_k$ and since $C_k \subset A_{k-1}$, we have that $|A_k|=|A_{k-1}|-|C_k|$ and the frequency $a_k$ is given by $a_k=\frac{|A_k|}{M}$. 

\begin{rem}
Let $n \in \N$ and $1 \leq k \leq n$. In order to define the frequencies $a_k$, we can interchange the action to construct sets $B_k$ and $C_k$ and frequencies $a_k$ do not change. This means that for any $k$, we can start from the previous set $A_{k-1}$ of $|A_{k-1}|=a_{k-1} M$ integers. We can apply the previous step (step k-1) to $A_{k-1}$ which gives the set $B_k$ satisfying $|B_k|=|A_{k-1}|a_{k-1}$. We can consider the subset $C_k$ of $B_k$ by keeping every number in $B_k$ from the first one by counting up in increments of $k$ and crossing out the remaining ones. The set $A_k$ is given by $A_{k-1} \setminus C_k$, we have that $|A_k|=|A_{k-1}|-|C_k|$ and $a_k=\frac{|A_k|}{M}$.
\end{rem}

\subsection{A good integer $M$.} Given $n \in \N$ and $1 \leq k \leq n$ it is necessary to determine good integers $M$ to make the construction above. 
\begin{defn}
A good integer $M$ is given by a positive integer such that in any step $k$ we can sieve the set $A_{k-1}$ correctly. That is, $|A_{k-1}|$ must be multiple of $k$ in order to construct $B_k$ and since we need to apply the previous step in $B_k$, we need that $|A_{k-1}||B_k|$ is multiple of $M$, so $\frac{|A_{k-1}|^2}{k M}$ must be an integer. 
\end{defn}

The following lemma is clear from definition of $A_k$.
\begin{lema} \label{lema exp}
Let $n \in \N$ and a good integer $M$ for $n$. For any $1 \leq k \leq n$ we have that
\begin{eqnarray} \label{expr}
|A_k|=|A_{k-1}|-\frac{|A_{k-1}|^2}{k M}.
\end{eqnarray}
\end{lema}

Fix $n \in \N$. We look for a good integer $M$ for $n$. Notice that the set $A_1$ satisfies that $|A_1|=1M=M$ so to assure that $|A_1|$ is an integer, it is sufficient to consider $M=1$. 

Now, $|A_2|=|A_1|-|A_1|^2/2M=M-M^2/M=M/2$, so it is sufficient to consider $M=2$ to assure that $|A_2|$ is an integer. %Analogously, $|A_3|=|A_2|-|A_2|^2/3M=M/2-M^2/12M=5M/12$, so it is sufficient to consider $M=12$ to assure that $|A_2|$ is an integer. 
Fix $2 < k \leq n$. If $|A_{k-1}|=\frac{\alpha}{\beta} M$, we have that $M=\beta$ is sufficient. Since
$$|A_k|=|A_{k-1}|-\frac{|A_{k-1}|^2}{k M}=\left(  \frac{\alpha}{\beta} -\frac{\alpha^2}{\beta^2 k}\right)M = \frac{\alpha \beta k- \alpha^2}{\beta^2 k}M,$$
\noindent it is clear that from step $k-1$ to step $k$ we change $\beta$ by $\beta^2 k$. Hence, we can give a sequence of good integers $M$ for all $n \in \N$.
\begin{cor}
A good integer for $n=1$ is given by $M_1=1$ and the recurrence $M_n=n \beta^2 =n M_{n-1}^2$ gives good integers $M_n$ for any $n \in \N$. 
\end{cor}

Notice that for any $n \in \N$, any multiple $M$ of $M_n$ is also a good integer for $n$ since we can still sieve correctly the set $\{1,2,\cdots,M\}$. 

From the recurrence, it is easy by induction that
\begin{prop}
Let $n \in \N$. Then,
$$M_n=2^{2^{n-2}} 3^{2^{n-3}} 4^{2^{n-4}} \cdots (n-1)^{2^{1}} n^{2^0}.$$
\end{prop} 

\subsection{The frequencies $a_n$.} Notice that the frequency $a_n$ is given by $\frac{|A_n|}{M}$. It is clear that this expression does not depends on $M$ since it is homogeneous of degree $0$ on $M$. We have used good integers $M$ in order to sieve the sets $A_n$ correctly but we can avoid this by considering upper limits when $M \rightarrow \infty$:
\begin{prop}
For any $n \in \N$ we have
$$a_n=\limsup_{M \rightarrow \infty} \frac{|A_n|}{M}.$$ 
\end{prop}

Dividing by $M$ in expression (\ref{expr}) in Lemma \ref{lema exp}, we can express $a_n$ by $a_1=1$ and
$$a_{n}=a_{n-1} \left(1-\frac{a_{n-1}}{n} \right).$$
\noindent Hence,
\begin{defn} \label{def pr}
For any $n \in \N$, we define the frequency $a_n$ by $a_1=1$ and
$$a_{n}=a_{n-1} \left(1-\frac{a_{n-1}}{n}\right)$$
\noindent for any $n \geq 2$. We will denote by $c(n)$ the sum $\sum_{k=1}^n a_k$.
\end{defn}

\begin{rem}
Notice that
$$a_{n}=a_{n-1}-\frac{a_{n-1}^2}{n},$$
\noindent so $a_{n}-a_{n-1}=-\frac{a_{n-1}^2}{n}$. Roughly speaking, this means that the ``derivative" of $``a_{n-1}"$ equals $``-a_{n-1}^2/n"$. Notice that the function $y=\frac{1}{\log x}$ also satisfies $y'=\frac{-1}{x}\frac{1}{\log^2 x}$ so, as we will show, the behaviour of $a_n$ will be similar to $\frac{1}{\log n}$. 
\end{rem}

\section{Results}
The following lemma is an easy calculation.

\begin{lema} \label{lema}
For any $n \geq 2$, we have that
$$\frac{1}{a_{n}}=\frac{1}{a_{n-1}}+\frac{1}{n-a_{n-1}}.$$
\end{lema}

\begin{lema} \label{lema2}
The series 
$$\sum_{k=2}^{\infty} \frac{a_{k-1}}{k(k-a_{k-1})}$$
\noindent is convergent to $S \approx 0.662834$.
\end{lema}

\dem It is clear that the series is convergent since $\sum_{k=1}^\infty \frac{1}{k^2}$ is convergent and $0 \leq a_k \leq 1$. An easy computation estimates the value of the sum. \qed \bigskip

The following proposition states that $1/a_n$ is very close to $\log n$. Indeed, 
\begin{prop} \label{prop1}
Let $b_n=\frac{1}{a_n}-\log n$ for any $n \in \N$. Then,
\begin{itemize}
\item[a)] The sequence $(b_n)$ is bounded:
$$\frac{1}{2}+\gamma < b_n < 1+\gamma.$$

\item[b)] The sequence $b_n$ is increasing.

\item[c)] The sequence $b_n$ is convergent to $\gamma + S \approx 1.24005$.
\end{itemize}
\end{prop}

\dem a) Using Lemma \ref{lema}, we have that for any $k \geq 2$,
\begin{eqnarray} \label{equ}
b_{k}-b_{k-1}=\frac{1}{k-a_{k-1}}+\log (k-1) - \log k.
\end{eqnarray}
\noindent Since $b_1=1$ and considering equality (\ref{equ}) for indexes from $2$ to $n$, we do a telescoping sum and obtain
\begin{eqnarray} \label{equ2}
b_n=1+\sum_{k=2}^{n} \frac{1}{k-a_{k-1}}-\log n.
\end{eqnarray}
Bearing in mind that $a_1=1$ and $0 \leq a_n \leq 1$ for any $n \geq 2$, we have making easy calculations that
$$\frac{1}{2}+H_{n}-\log n < b_n < 1-\frac{1}{n}+H_{n}-\log n,$$
\noindent so bearing in mind inequality (\ref{desig2}), we have 
$$\frac{1}{2} +\gamma < b_n < 1+\gamma$$ 
\noindent and we are done. \medskip

b) We will show that $b_{n} > b_{n-1}$ for any $n \geq 2$. This is true if and only if
$$\frac{1}{a_{n}}-\log n-\left( \frac{1}{a_{n-1}}-\log (n-1) \right) >0,$$ 
\noindent if and only if
$$\frac{1}{a_{n}}-\frac{1}{a_{n-1}} -\log \left(\frac{n}{n-1} \right) >0.$$
\noindent Using definition of the sequence $a_n$, this inequality is equivalent to 
$$\frac{a_{n-1}}{n a_{n}}-\log \left(\frac{n}{n-1} \right) >0 \ \mbox{ if and only if }
 \ a_{n-1} > \log \left(\frac{n}{n-1} \right)^n a_{n}.$$
\noindent Notice that $(a_n)$ is decreasing since $a_{n}-a_{n-1}=-\frac{a_{n-1}^2}{n} <0$, so
$$a_{n-1} > a_{n} > \log \left(\frac{n}{n-1} \right)^n a_{n},$$
\noindent where last inequality is true since $\left(\frac{n}{n-1} \right)^n$ is an increasing sequence which tends to $e$ and $\log x$ is an increasing function on its domain. \medskip

c) Since $(b_n)$ is increasing and bounded, it is convergent to some limit $\ell$. It is clear that 
$$b_n=H_{n} - \log n + \sum_{k=2}^{n} \left( \frac{1}{k-a_{k-1}} - \frac{1}{k} \right)=H_{n} - \log n + \sum_{k=2}^{n} \frac{a_{k-1}}{k(k-a_{k-1})}.$$
\noindent Last sum is a convergent series to $S$ by Lemma \ref{lema2}. Then,
$$\lim_{n \rightarrow \infty} b_n=\lim_{n \rightarrow \infty} H_{n}- \log n+\sum_{k=2}^{n} \frac{a_k}{k(k-a_k)}=\gamma+S$$
\noindent and we are done. \qed \bigskip

\begin{cor}
We have that $\frac{1}{a_n} \sim \log n$.
\end{cor}

%\dem It is clear since, dividing by $\log n$, we have that 
%$$\frac{\frac{1}{2}-\frac{1}{n}}{\log n} < \frac{\frac{1}{a_n}}{\log n }-1 < \frac{1+\gamma}{\log n},$$
%\noindent so taking limits when $n \rightarrow \infty$ we are done. \qed \bigskip

Notice that from equality (\ref{equ2}) in Proposition \ref{prop1}, we have that
\begin{cor} \label{corol}
For any $n \geq 2$ we have
$$\frac{1}{a_n}=1+\sum_{k=2}^{n} \frac{1}{k-a_{k-1}}.$$
\end{cor}

The following lemma is an easy consequence of $\pi(x) \sim Li(x)$ and integral calculus.
\begin{lema}
Consider the function 
$$q(n)=\sum_{k=2}^{n} \frac{1}{\log k}$$
\noindent for any $n \in \N$. Then, $\pi(n) \sim q(n)$.
\end{lema}

\begin{lema} \label{series}
For any real number $C$ and an integer number $m \geq 2$, we consider the sequence 
$$x_{m,C}=\sum_{n=2}^{m} \frac{1}{\log(n)+C}.$$ 
\noindent Then, $x_{m,C}$ is divergent and $x_{m,C_1} \sim x_{m,C_2}$ for any $C_1, C_2 \geq 0$. 
\end{lema}

\dem It is well-known that $x_{m,1}$ is a divergent series and by the limit criterium, $x_{m,C}$ is also divergent for any $C \geq 0$. Hence, the limit
$$\ell = \lim_{m \rightarrow \infty} \frac{x_{m,C_1}}{x_{m,C_2}}$$
\noindent can be calculated using the Stolz criterium, so we obtain that 
$$\ell= \lim_{m \rightarrow \infty} \frac{\frac{1}{\log(m+1)+C_1}}{\frac{1}{\log(m+1)+C_2}}=\lim_{m \rightarrow \infty} \frac{\log(m+1)+C_2}{\log(m+1)+C_1}=1.$$

\begin{teo}
The sum $c(n)=\sum_{k=1}^n a_k$ satisfies $c(n) \sim \pi(n)$.
\end{teo}

\dem By Corollary \ref{corol},
$$\frac{1}{2} + H_{k} < \frac{1}{a_k} < 1+ H_k-\frac{1}{k}.$$
\noindent By inequality \ref{desig2}, we have that $\log k + \gamma < H_k < \log k + \gamma + \frac{1}{k},$ so 
$$1+ \log k < \frac{1}{2}+ \log k + \gamma \leq \frac{1}{a_k} \leq 1+\log k + \gamma < 2+\log k$$
\noindent for any $k \geq 2$. Hence, we have that
$$\frac{1}{\log k+2} \leq a_k \leq \frac{1}{\log k+1},$$
\noindent so bearing in mind that $a_1=1$, we have
%$$\sum_{k=1}^n \frac{1}{\log k+2} \leq \sum_{k=1}^n a_k \leq \sum_{k=1}^n \frac{1}{\log k+1},$$
%\noindent which can be written as
\begin{eqnarray} \label{desig}
\sum_{k=1}^{n} \frac{1}{\log k+2} \leq \sum_{k=1}^n a_k \leq \sum_{k=1}^{n} \frac{1}{\log k+1}.
\end{eqnarray}
 
Using Lemma \ref{series},
dividing by $x_{n,0}$ in inequality (\ref{desig}), we obtain that 
$$1 \leq \lim_{n \rightarrow \infty} \frac{c(n)}{\sum_{k=2}^n \frac{1}{\log k}} \leq 1,$$
\noindent and by the Sandwich criterium, we are done. \qed \bigskip

\begin{cor}
We have $\frac{n}{\pi(n)} \sim \frac{1}{a_n}$ and $\pi(n) \sim n a_n$.
\end{cor}

\dem By Proposition \ref{prop1}, we have that $\frac{1}{a_n} \sim \log n$. Since $\frac{n}{\log n} \sim \pi(n)$, we have that $\frac{n}{\pi(n)} \sim \log n \sim \frac{1}{a_n}$. The other statement is also clear. \qed \bigskip

\subsection{Evaluating $n/\pi(n)$.} It was first Chebyshev who proved (see \cite{C}) that there exists $x_0 \in \N$, $c_1 \approx 0.92$ and $c_2 \approx 1.1$ such that 
$$c_1 \frac{x}{\log x} \leq \pi (x) \leq c_2 \frac{x}{\log x}$$
\noindent for any $x \geq x_0$. The following result was obtained in \cite{RS}:
$$\frac{x}{\log x-1/2} < \pi(x) < \frac{x}{\log x-3/2},$$
\noindent where first inequality is true for $x \geq 67$ and the second one is true for $x \geq e^{3/2}$. L. Panaitopol (see \cite{P}) improved this result:
\begin{eqnarray} \label{pan}
\frac{x}{\log x-1+\frac{1}{\sqrt{\log x}}} < \pi(x) < \frac{x}{\log x-1-\frac{1}{\sqrt{\log x}}},
\end{eqnarray}
\noindent where first inequality is true for $x \geq 59$ and the second one for $x \geq 6$. 

Hence, we obtain the following result:
\begin{prop}
We have the following estimates for $n \geq 59$:
\begin{itemize}
\item[1)] We have that 
$$\left|\frac{n}{\pi(n)}-\log n +1 \right| \leq \frac{1}{\sqrt{\log n}},$$
\noindent so 
$$\lim_{n \rightarrow \infty} \left| \frac{n}{\pi(n)}-\log n+1  \right| =0.$$

\item[2)] We have that 
$$\left|\frac{1}{a_n}-\frac{n}{\pi(n)} \right| \leq 2+ \gamma +\frac{1}{\sqrt{\log n}} < 3.$$

\item[3)] We have that
$$\lim_{n \rightarrow \infty} \frac{1}{a_n}-\frac{n}{\pi(n)}=\gamma+S+1 \approx 2.24005.$$
\end{itemize}
\end{prop}

\dem The first one is true by inequality \ref{pan} from the Panaitopol result. The second result is an easy consecuence from the first one and Proposition \ref{prop1}. To prove the third one, notice that from the first result we have that
$\log n-\frac{n}{\pi(n)}-1 \rightarrow 0$ when $n \rightarrow \infty$ so
$$\frac{1}{a_n}-\frac{n}{\pi(n)} = \frac{1}{a_n}-\log n + 1+\log n -\frac{n}{\pi(n)}-1$$
\noindent which tends to $\gamma+S+1$ when $n \rightarrow \infty$ by Proposition \ref{prop1}. \qed \bigskip

%\subsection{Remarks and Questions.}
%
%%We are interested in estimating the differences 
%%$$d_1=\pi(n)-s(n),$$
%%$$d_2=\frac{1}{a_n}-\log n,$$
%%$$d_3=\frac{n}{\pi(n)}-\frac{1}{a_n}.$$
%%
%%Numerical calculations suggests that $d_3 \in [-2.32,-2.29]$ for $n$ between $10^4$ and $10^{10}$.
%Maybe it is interesting for something:
%\begin{itemize}
%%\item[1)]  . \medskip
%
%\item[2)] Notice that 
%$$\frac{1}{1-\frac{a_n}{n}}=\sum_{k=0}^\infty \left( \frac{a_n}{n} \right)^k,$$
%\noindent so from the definition of the recurrence we have
%$$\frac{1}{a_n} \sum_{k=0}^\infty \left( \frac{a_n}{n} \right)^k = \frac{1}{a_{n+1}}.$$
%\end{itemize}

\bibliographystyle{amsplain}

\end{document}